# CORRECTION

## STRONG INVARIANCE PRINCIPLES FOR SEQUENTIAL BAHADUR–KIEFER AND VERVAAT ERROR PROCESSES OF LONG-RANGE DEPENDENT SEQUENCES

*The Annals of Statistics* (2006) **34** 1013–1044

BY MIKLÓS CSÖRGŐ, BARBARA SZYSZKOWICZ AND LIHONG WANG

*Carleton University, Carleton University and Nanjing University*

Rafał Kulik of University of Sydney and Wrocław University has brought to our attention that Assumption A does not help in the proof of our Proposition 2.2, for what is really used in our method of proof is, actually, the uniform boundedness of $J_\tau(y)$ and that of its derivatives. However, if $X_i = \eta_i$ (cf. Remark 2.1), then $J_1(y) = -\phi(\Phi^{-1}(y))$, and thus, we have $J_1'(y) = -\Phi^{-1}(y)$, $J_1''(y) = -(\phi(\Phi^{-1}(y)))^{-1}$, and both are unbounded functions over the unit interval. One arrives at a similar conclusion in the case of $X_i = \eta_i^2$, that is, when $G(x) = x^2$ (cf. Remark 1.1) and $\tau = 2$. Consequently, the proof of Proposition 2.2 is not valid, unless we restrict ourselves to:

(R) intervals $y \in [a,b]$, $0 < a < b < 1$ instead of $y \in [0,1]$,
or assume that $F$ has finite support.

Hence, we conclude, for further use as well, the following observation.

REMARK. Instead of Assumption A being assumed in Proposition 2.2, for the validity of its present proof, we must assume the above restriction (R).

This Remark now automatically applies also to Propositions 2.3, 2.4 and 2.5, as well as to Theorems 2.2, 2.3 and 2.4.

We note in passing that Theorem 2.1 continues to hold true as stated, that is, under the assumptions of Corollary 2.1.

As a consequence of our comments so far, and due to the definition of the sequential uniform Vervaat error process $V_n(\cdot,\cdot)$ as in (1.10), we conclude that Theorem 3.1, as well as Proposition 3.1, continue to hold true, provided that $F$ has finite support. The same holds true for Theorems 3.2 and 3.3, in which the constant $2^{5/2}$ should be replaced with $2^{3/2}$. The reason for this is that there is a mistake in Proposition 3.2, as stated. The correct version is as follows.

PROPOSITION 3.2. *Under the assumptions of Proposition 2.2, and assuming that $F$ has finite support, as $n \to \infty$, we have*

$$\sup_{0 \le t \le 1} \sup_{0 \le s \le 1} \left| Z_n(s,t) + d_n^{-2}[nt]^{-1}(V(s,nt))^2 J'_\tau(s) \sum_{i=1}^{[nt]} H_\tau(\eta_i)/\tau! \right|$$
$$= O(n^{-\tau D} L^\tau(n)(\log \log n)^{2\tau}) \qquad a.s.$$

PROOF. By the same argument as in the end of the proof of Proposition 2.1, it suffices to show that the above bound is valid for $t = 1$. We have

$$V(s - wn^{-1}V(s,n), n) - V(s,n)$$
$$= -V'(s,n)wn^{-1}V(s,n) + O(V''(s,n)(wn^{-1}V(s,n))^2) \qquad \text{a.s.}$$

Thus, bearing in mind the definition of $Z_n(\cdot, \cdot)$ [cf. (3.1)],

$$\sup_{0 \le s \le 1} \left| Z_n(s,1) + 2d_n^{-2}n^{-1}(V(s,n))^2 J'_\tau(s) \sum_{i=1}^n H_\tau(\eta_i)/\tau! \int_0^1 w\,dw \right|$$
$$= O(d_n^{-2}n^{-2}V''(s,n)(V(s,n))^3) = O\left(d_n^{-2}n^{-2}\left(\sum_{i=1}^n H_\tau(\eta_i)/\tau!\right)^4\right)$$
$$= O(n^{-\tau D} L^\tau(n)(\log \log n)^{2\tau}) \qquad \text{a.s.,}$$

where, on assuming that $F$ has finite support, we have made use of $\sup_{0 \le s \le 1} J''_\tau(s) < \infty$ and the law of iterated logarithm for partial sums [cf. (2.7)]. This completes the proof. $\square$

Concerning Section 4, Proposition 4.2 remains valid on assuming the condition (R), in addition to (i)–(iv) with (v) or (v′) for $F$. However, its almost sure bounds in (4.5) hold true only for $\gamma \in (0, 1]$. If $\gamma > 1$, the indicated bounds are valid only in probability. The reason for this is as follows. In [1] the following property of uniform order statistics was used in the i.i.d. case: $\liminf_{n \to \infty} n(\log n)^{1+\varepsilon} U_{1,n} = \infty$ almost surely for all $\varepsilon > 0$. There is no such result available in the long-range dependent case.

We note in passing that the approximation of the general quantile process by the uniform quantile process as stated in (4.7) remains valid under the conditions (i)–(iii) on $F$ without assuming also (R).

The results for the general Bahadur–Kiefer process (Theorems 4.1 and 4.2) cannot be concluded from those for the uniform Bahadur–Kiefer process, as claimed in Remark 4.1 of the paper in view of (4.12). The reason for this is that the rate at which the uniform Bahadur–Kiefer process behaves (cf. Theorem 2.3) is, mutatis mutandis, the same as that of approximating the general quantile process $\rho_n(y,t)$ by the uniform quantile process [cf. (4.7) in combination with (4.9) and (4.10)]. Moreover, in view of results in [2], we conjecture that the limiting process in Theorem 4.1 is to be changed to $1/2$ times that of the claimed process, that is, to

$$\frac{1}{(2-\tau D)(1-\tau D)} J_\tau(y) J'_\tau(y) Y_\tau^2(t).$$

The corresponding conjecture also applies to the remaining results of Theorems 4.1 and 4.2. In view of this comment and the scaling constants in the invariance principle for the uniform Bahadur–Kiefer process (cf. Theorem 2.3), we conclude that the big "$O$" almost sure rates in Proposition 4.2 cannot be replaced with "$o$." In this regard, we note that Wu in [2] considered related approximations in case of nonsubordinated linear processes on $[a,b]$, $0 < a < b < 1$.

**Acknowledgments.** The authors wish to thank Rafał Kulik for bringing these points to their attention, for helping to correct their oversights and for reviewing their results in this regard.

M. Csörgő
B. Szyszkowicz
School of Mathematics and Statistics
Carleton University
1125 Colonel By Drive
Ottawa, Ontario
Canada K1S 5B6
E-mail: mcsorgo@math.carleton.ca
　　　　bszyszko@math.carleton.ca

L. Wang
Department of Mathematics
Nanjing University
Nanjing 210093
China
E-mail: lhwang@nju.edu.cn